\documentclass[a4paper,11pt]{amsart}
\usepackage{amssymb}
\usepackage{latexsym}
\usepackage{amsthm}
\usepackage{hyperref}
\newtheorem{thm}{Theorem}
\newtheorem{cor}[thm]{Corollary}
\newtheorem{lem}[thm]{Lemma}

\theoremstyle{remark}

\newtheorem{exmp}[thm]{Example}

\theoremstyle{definition}

\newenvironment{pf}{\par\noindent{\bf Proof.}\enspace\ignorespaces}{\qed\par\par}

\def\qed{\hfill $\Box$}

 \newcommand{\bZ}{{\mathbb{Z}}}
 \newcommand{\bP}{{\mathbb{P}}}
 
 \newcommand{\bF}{{\mathbb{F}}}
 
 \newcommand{\cP}{{\mathcal{P}}}

 \DeclareMathOperator{\Jac}{Jac}
 \DeclareMathOperator{\Div}{Div}
 
 \DeclareMathOperator{\PGL}{PGL}
 
 \DeclareMathOperator{\AS}{AS}
 \DeclareMathOperator{\Hyp}{Hyp}

 \DeclareMathOperator{\Trig}{Trig}

\begin{document}

\title{A census of all genus 4 curves over the field with 2 elements}

\author[Xavier Xarles]{Xavier Xarles }
\address{Departament de Matem\`atiques\\Universitat Aut\`onoma de
Barcelona\\08193 Bellaterra, Barcelona, Catalonia}
\email{xarles@mat.uab.cat}

\begin{abstract} We explain how we computed equations for all genus 4 curves defined of the field $\bF_2$ up-to-isomorphism, and some of the data we obtained. We give descriptions also of nice models for genus 4 curves over characteristic 2 fields, in both the hyperelliptic case and the trigonal case. 
\end{abstract}

\maketitle

In this short note we explain how we computed equations for all genus 4 curves defined of the field $\bF_2$ with 2 elements up-to-isomorphism. 
We also give tables of some of the results as well as some interesting examples of curves. 

Recall that, any curve of genus 4 over a field either it is hyperelliptic or it
is trigonal. In ths first case it has an hyperelliptic equation, and in the second can be given as intersection of a quadric and a cubic in the projective space $\bP^3$. We computed both cases separately.

\subsection*{Acknowledgements}

The author is supported by MTM2016-75980-P. I thank Joan Carles Lario for proposing the problem that let to this project, and for some detailed dicussions, and Joaquim Ro\'e, Francesc Bars and Marc Masdeu for some discussions related to the problem. 

\subsection*{Notations}

For a given (finite) field $K$, we denote by $\Hyp_g(K)$ the set of
all hyperelliptic curves of genus $g$ defined over $K$ modulo
isomorphism over $K$,  and by $\Trig_g(K)$  the set of
all  trigonal curves of genus $g$ defined over $K$ modulo
isomorphism over $K$ . Our aim is to describe an algorithm to compute
$\Hyp_4(\bF_2)$ and $\Trig_4(\bF_2)$, and compute them explicitly. 

We denote by $K[x]_n$ the set of all  all nonzero polynomials with degree $\le n$. 
 Note that trivially  $\bF_2[x]_n$ consists on the set of monic polynomials (since monic implies non-zero). 

For any characteristic 2 field $K$,
we denote by $$\AS(K):=\{a^2+a \in K\ | \ a\in K\}\subset K$$ the
Artin-Schreier subgroup. Observe that $K/\AS(K)$ classifies all
degree 2 extensions of the field; if $L$ is a finite field, then $K/\AS(K)$ has only two elements.

\section{How we computed all isomorphism classes: The hyperelliptic
case}

We recall some results on hyperelliptic curves over characteristic 2
fields. For another but equivalent approach to classify such hyperelliptic curves see \cite{NaSa}. 

Recall that (see for example Proposition 7.4.24 in \cite{Liu}) any
hyperelliptic curve has a model of the form
$$ y^2+q(x)y=p(x)$$
with $p$ and $q$ polynomials verifying that $$2g+1\le
\max\{2\deg(q),\deg(p)\}\le 2g+2.$$ 

Doing the change of variables $y\mapsto \alpha y$, we could also suppose
that the polynomial $q(x)$ is monic. In case $K$ has characteristic 2, we must have $q(x)\ne 0$ 
(in order to have an irreducible curve). We will call such a model a standard model. 

Consider the set $K[x]_n$ of all nonzero polynomials
with degree $\le n$ over $K$. Given any matrix $$A=\begin{pmatrix} a & b \\
c & d\end{pmatrix}\in \PGL(2,K),$$ and any polynomial $q(x)\in
 K[x]_n$, we denote by $$\psi_n(A)(q):=(cx+d)^n q\left(
\frac{ax+b}{cx+d} \right).$$ Then, $\psi_n$ induces an action of
$\PGL(2,K)$ into the set  $K[x]_n$. We denote by
$\overline{K[x]_n}= K[x]_n/\PGL(2,K)$ the quotient set.

\begin{lem}\label{qiso} Let $H_1$ and $H_2$ be two hyperelliptic curves over $\bF_2$ 
given by an equation is standard form
$$ y^2+q_i(x)y=p_i(x)$$
 with $2g+1\le \max\{2\deg(q_i),\deg(p_i)\}\le 2g+2$
and $q_i(x)$ monic, for $i=1,2$. If $H_1\cong H_2$, then there
exists a matrix $A\in \PGL(2,\bF_2)$ such that
$q_2(x)=\psi_{g+1}(A)(q_1(x)).$ 
\end{lem}
\begin{pf}
We have that $H_1$ and $H_2$ must be related by a change of variable of the form 
$$(x,y) \mapsto \left(\frac{ax+b}{cx+d},\frac{r(x)+y}{(cx+d)^{g+1}}\right), \text{ where } A=\begin{pmatrix} a & b \\
c & d\end{pmatrix}\in \PGL(2,\bF_2),$$
and $r(t)$ is a polynomial of degree $\le g+1$  (see 7.4.33 in \cite{Liu}). Denote by $\varphi(r,A)$ such a change of variables. Note that 
$ \varphi(r,A)= \varphi(0,A) \varphi(r,1)$.

But in characteristic $2$ one has that 
$$(r(x)+y)^2=r(x)^2+y^2$$
so the coefficient of $y$ in the hyperelliptic equation after a base change of the type $\varphi(r,1)$ is invariant, and only the action of the linear group can affect it. A simple computation shows the result. 
\end{pf}

In other words, the lemma says that the map
$$\Hyp_g(\bF_2) \to\overline{\bF_2[x]_{g+1}}$$
defined by sending any hyperelliptic curve in standard form 
 $y^2+q(x)y=p(x)$ with $p$ and $q$ polynomials
verifying that $2g+1\le \max\{2\deg(q),\deg(p)\}\le 2g+2$ to the
class of $q$ is well defined. One can also see that it is
surjective in general (a fact we will see in practice in our case). 

The following lemma determines a uniquely determined representative
for each class in $\overline{\bF_2[x]_{g+1}}$ for $g=4$.

\begin{lem}\label{qs} The following list $Q_4(\bF_2)$ contains a complete set of representatives
of $\overline{ \bF_2[x]_{5}}$:  with degree $\le 2$,
$$    1,\;
    x,\; 
    x^2,\;
    x(x+1)=x^2 + x,\; 
    x^2 + x + 1,$$
with degree 3
$$
 x^2(x+1)=x^3+x^2,\;
x(x^2+x+1)=x^3 +x^2+x,\;
    x^3 + x + 1 ,$$
with degree 4
$$ x(x+1)(x^2+x+1)=x^4 + x,\;
(x^2+x+1)^2= x^4 + x^2 + 1,\; x(x^3+x+1)=x^4 + x^2 + x,$$
$$    x^4 +
x + 1,\;
    x^4 + x^3 + 1,$$
and with degree 5
$$  x^5 + x + 1=(x^3+x^2+1)(x^2+x+1) \text{ and }
    x^5 + x^2 + 1.
$$
\end{lem}

\begin{pf}
Consider the divisor $D_q:=Z(q(x))+(g+1-\deg(q)) \infty$ in $\bP^1_{\bF_2}$,
where $Z(q)$ denotes the zero divisor of $q$; it gives an element of $\Div_{g+1}( \bF_2)$, the divisors of degree $g+1$ 
by points in $\bP^1_{\bF_2}$ . One easily shows that
the action of $\PGL(2, \bF_2)$ in $\overline{ \bF_2[x]_{g+1}}$ translates to the natural
action on the divisors group on $D_q$. This gives an identification
of $\overline{ \bF_2[x]_{g+1}}$ with $\Div_{g+1}( \bF_2)/\PGL(2, \bF_2)$ form by the
degree $g+1$ divisors modulo the action of $\PGL(2, \bF_2)$. 

It is well
known that there exists a unique matrix in $\PGL(2, \bF_2)$ sending three
given distinct $K$-rational points $p_0$, $p_1$ and $p_\infty$ to
$0$, $1$ and $\infty$, respectively. So, any divisor containing only
three rational points can be written as $n_\infty \infty + n_0 0 +
n_1 1$ with $n_\infty\ge n_0 \ge n_1$. Since $\bP^1(\bF_2)$ contains
only this three points, this accounts for all divisors form by
$K$-rational points. The possibilities when $g=4$ are
$$(n_\infty,n_0,n_1)=(5,0,0), (4,1,0), (3,2,0),(3,1,1) \text{ or } (2,2,1).$$
This accounts for all the cases that the polynomial $q(x)$ factorizes completely over $\bF_2$.

It remains to consider the cases when the divisor contains a degree
d point, for some $2\le d\le 5$. If
it has a degree 2 point $\alpha_2$, and the rest $K$-rational
points, modulo the action of $\PGL(2,\bF_2)$ we have that
$$D_q=3\infty + (\alpha_2), 2\infty+ 0 +(\alpha_2), \infty+ 0 +1
+(\alpha_2)  \text{ or } \infty + 2(\alpha_2).$$ If it has a degree 3 point
$\alpha_3$, we get
$$D_q=2\infty + (\alpha_3), \infty+ 0 +(\alpha_3)  \text{ or }
(\alpha_2)+(\alpha_3).$$ Note that, although there are two possible
elections of a degree 3 point (given as a root of $x^3 + x + 1$ or a
root of $x^3+x^2+1$), both can be interchanged by either $x\mapsto
1/x$ or $x\mapsto x+1$, which preserve the rest of the divisor
$D_q$.

If it has a degree 4 point $\alpha_4$, only one option is possible,
$D_q=\infty + (\alpha_4)$. Here we have that the $3$ distinct
options for $\alpha_4$, although they can be interchanged by the
action of $\PGL(2,\bF_2)$, they do not give equivalent divisors
$D_q$. For example the ones given by $x^4+x+1$ and $x^4+x^3+1$, as
they can only be interchanged by $x\mapsto 1/x$ (using elements of
$\PGL(2,\bF_2)$), they change the divisor $D_q$. So two distinct
divisors appear.

Finally, with a degree $5$ point $\alpha_5$, the only possibility is
$D_q=(\alpha_5)$. In this case all the six distinct points can be
interchanged by $\PGL(2,\bF_2)$, and so only one $D_q$ remains.
\end{pf}

Now, for any given $q \in Q_4(\bF_2)$, we will denote by $G_q$ the
stabilizer of $q$ with respect to the action of $\PGL(2,\bF_2)$.

\begin{exmp} For the polynomials 
$$q= x, x^2, x^3 + x^2+x, 
x^3 + x + 1,
x^4 + x^2 + x,
x^4 + x^3 + 1,
x^5 + x^2 + 1$$
the group $G_q$ is trivial. 

On the other hand, for $q=x^4 + x$, we have $G_q=PGL(2,\bF_2)$, and for $q=x^5 + x + 1$ it is the only subgroup of order 3. All the other case it is a subgroup of order 2. 
\end{exmp}

Our aim now is to compute the set $\cP_q$ of polynomials $p(x)$ with
degree $\le 2g+2=10$ such that the curve given be the equation
$y^2+q(x)y=p(x)$ is of genus $g=4 $, modulo isomorphism. I.e., we want to
compute the fiber of $q$ with respect to the map $$\Hyp_4(\bF_2) \to
\overline{\bF_2[x]_{5}}\cong Q_4(\bF_2).$$ In order to do it
systematically we need first to know when the curve is non-singular
and it has genus $4$, and then to know when two such curves are
isomorphic; in the following lemma we explicitly state the criterion to know the 
isomorphism class of a curve of the form $y^2+qy=p$ for a fixed $q$. The following result 
it is easily deduced from the explicit description of the isomorphisms we recalled in the proof of lemma \ref{qiso}.

\begin{lem}\label{isoq} Let $H_1$ and $H_2$ be two hyperelliptic curves of genus 4 over $\bF_2$ 
given by equations in standard form
$$ y^2+q(x)y=p_i(x)$$
(with $9\le \max\{2\deg(q_i),\deg(p_i)\}\le 10$)
for $i=1,2$. If $H_1\cong H_2$, then there
exists a matrix $A \in G_q$ and a polynomial $r(x)\in \bF_2[x]$ with degree $\le g+1=5$ such that 
$p_2(x)=\psi_9(A)(p_1(x)+r(x)^2+q(x)r(x))$.
\end{lem}

The first condition we need is on the degree to verify that
 $$9\le \max\{2\deg(q),\deg(p)\}\le 10.$$
This condition says that, if  $\deg(q)<5$, then $\deg(p)=9$ or $10$,
while if $\deg(q)=5$ it says nothing.

\begin{lem}\label{affine} Fix $g\ge 1$ and a field $K$ of characteristic $=2$. Given $q(x)\in K[x]$ 
monic of degree $\le g+1$ and non-constant 
and $p(x)\in K[x]$ 
with degree $\le 2g+2$ with  $2g+1\le \max\{2\deg(q),\deg(p)\}\le 2g+2$. 
Then the hyperelliptic equation  $y^2+q(x)y=p(x)$ gives a curve of genus $g$ 
if and only if $$\gcd\left(q(x),p'(x)^2+q'(x)^2p(x)\right)=1$$
and 
either $\deg(q)=g+1$ or $a_{2g+1}^2\ne a_{2g+2}b_g^2$, where $p(x)=\sum_{i=0}^{2g+2} a_ix^i$ and 
$q(x)=\sum_{i=0}^{g+1} b_ix^i$. 
\end{lem}
\begin{pf}
The first condition on the gcd's concerns the non-singularity on the affine points. A point $(a,b)\in L^2$ is 
singular for the equation $f(x,y)=0$ if and only if 
$$ f(a,b)=0, \ \frac{\partial f}{\partial x}(a,b)=0 \ \text{ and } \frac{\partial f}{\partial y}(a,b)=0$$
In our case, considering $f(x,y)=y^2+q(x)y+p(x)$, we get the condition 
$$q'(a)b=p'(a), \ q(a)=0 \text{ and } b^2+q(a)b=p(a)$$
From this we get that $b^2=p(a)$, and hence $p'(a)^2=q'(a)^2p(a)$. Hence $a$ has to be a simultaneous root of $q$ and of $p'^2+q'^2p$. Reversely, for such a root $a$, either $q'(a)\ne 0$, and taking $b=p'(a)(q'(a)^{-1})$ we  get a singular point $(a,b)$, or $q'(a)=0$, hence $p'(a)=0$, and any root $b$ of $y^2=p(a)$ gives a singular point $(a,b)$. 

The second conditions appear to insure that the points at infinity are non-singular. If we do a base change corresponding to $(x,y)\mapsto (1/x,y/x^{g+1})$, the points at infinity correspond to the points with $x=0$ of the new equation $y^2+\bar q(x)y+\bar p(x)=0$. Now, the conditions to have $(0,b)$ as a singular point is $\bar q(0)$, and $\bar p'(0)^2+\bar q'(0)^2\bar p(0)=0$,  which translates in to the original equation as to $b_{g+1}=0$ (so $\deg(q)<g+1$) and  $a_{2g+1}^2+ b_g^2a_{2g+2}=0$. 
\end{pf}

\begin{cor}\label{infinity} Under the hypothesis of the previous lemma, if moreover $\deg(p)\le 2g+1$, the points at infinity are always non singular.  \end{cor}
\begin{pf} If $\deg(p)<2g+1$, then necessarily $\deg(q)=g+1$, hence the points at infinity are non-singular. If $\deg(p)=2g+1$, then $a_{2g+1}\ne 0$, and $a_{2g+2}=0$, hence $a_{2g+1}^2\ne a_{2g+2}b_g^2=0$. 
\end{pf}

\begin{cor}\label{infinityg4} Consider $q\in Q_4(\bF_2)$ and a polynomial $p(x)=\sum_{i=0}^{2g+2} a_ix^i \in \bF_2[x]$. Then  the curve given by $y^2+q(x)y=p(x)$ is non-singular at infinity if 
\begin{enumerate}
\item If $\deg(q)\le 3$ and $a_{9}\ne 0$. 
\item If $\deg(q)=4$ and $a_{10}\ne a_{9}$. 
\item If $\deg(q)=5$. 
\end{enumerate}
\end{cor}

In the following corollary we collected explicit criterion for non-singularity for some of the polynomials in $ Q_4(\bF_2)$.
I will do some of the cases, leaving the others to the reader. 

\begin{cor} The conditions for a polynomial $p(x)=\sum_{i=0}^{2g+2} a_ix^i \in \bF_2[x]$ to define a genus $4$ curve $y^2+q(x)y=p(x)$ are
\begin{enumerate}
\item if $ q=1$:  $a_{9}=1$. 
\item  if  $ q=x$: $a_{9}=1$ and $a_1\ne a_0$.
\item  if  $ q=x^2$: $a_{9}=1$ and $a_1=1$.
\item if  $ q= x(x+1)$: $a_{9}=1$, $a_1\ne a_0$ and  $(p'+p)(1)=\sum_{i=0}^{5} a_{2i}=1$.
\item  if $ q= x^2(x+1)$: $a_{9}=1$,  $a_1=1$  and $(p'+p)(1)=\sum_{i=0}^{5} a_{2i}=1$.
\item  if $ q= x^2 + x + 1$: $a_{9}=1$, $(\sum_{i=0}^{2}a_i)+ (\sum_{i=5}^{8} a_i)=1$ and $(\sum_{i=1}^{5}a_{2i})=1$.
\end{enumerate}
\end{cor}

\begin{pf} The condition $a_9\ne 0$ is directly from Corollary \ref{infinityg4}. So we only need to insure that $\gcd\left(q(x),p'(x)^2+q'(x)^2p(x)\right)=1$ for any case. 

For example, in the case (1), the gcd condition is automatic. In case (2) it says $\gcd\left(x,p'(x)^2+p(x)\right)=1$, so $p'(0)\ne p(0)$, so $a_1\ne a_0$. And in case (3) that $\gcd\left(x,p'(x)^2\right)=1$, so $p'(0)\ne 0$, thus  $a_1=1$. 

In case (3) the condition on the gcd becames $\gcd\left(x(x+1),p'(x)^2+p(x)\right)=1$, hence $p'(0)\ne p(0)$ and $p'(1)+ p(1)\ne 0$. Similarly, in case (4) we get  $\gcd\left(x^2(x+1),p'(x)^2+x^2p(x)\right)=1$, hence $p'(0)\ne 0$ and again $p'(1)+ p(1)\ne 0$. 

Finally, the last case we get the condition $\gcd( x^2 + x + 1,p'^2+p)=1$. Denote $\alpha\in \bF_4$ such that $\alpha^2=\alpha+1$. Then the condition on the gcd becames $p'(\alpha)^2+p(\alpha)\ne 0$. We have 
$$p(\alpha)=\sum_{i=0}^{10} a_i\alpha^i=(a_0+a_3+a_6+a_9)+(a_1+a_4+a_7+a_{10})\alpha+(a_2+a_5+a_8)\alpha^2=$$
$$=(a_0+a_2+a_3+a_5+a_6+a_8+a_9)+(a_1+a_2+a_4+a_5+a_7+a_8+a_{10})\alpha$$
and 
$$p'(\alpha)^2=\sum_{i=0}^{4} a_{2i+1}\alpha^{4i+2}=(a_1+a_3+a_7+a_9)+(a_1+a_5+a_7)\alpha.$$
The condition $p'(\alpha)^2\ne p(\alpha)$ becames then 
$$
a_0+a_2+a_5+a_6+a_8\ne a_1+a_7 \text{ and } 
a_2+a_4+a_8+a_{10}\ne 0.$$
\end{pf}

We implemented the following method to list all hyperelliptic
genus 4 curves over a finite field $k$. For a fixed polynomial $q$ as
in the lemma \ref{qs}
\begin{enumerate}
\item We computed the group $G_q$.
\item Take $D_{10}:=\bF_2[x]_{10}$ the set of polynomials  with degree $\le 10$.
\item We choose a polynomial $p$ in $D_{10}$ and we compute the set $S_p$  obtained as $$S_p':=\{p(x)+r(x)^2+q(x)r(x)\mid r(x)\in \bF_2[x]_5\},$$ and $$S_p:=\left\{ (cx+d)^{10}f( \frac{ax+b}{cx+d})\mid \begin{pmatrix} a&b\\ c&d \end{pmatrix} \in G_q  \text{ and } f\in S_p'\right\}.$$
\item We choose an element $f\in S_p$, we add to the set $V_q$ and we remove $S_p$ from $D_{10}$. We repeat from step (2) until $D_{10}$ is empty.
\item For each $f\in V_q$, we determine if they equation $y^2+qy=p$  gives an hyperelliptic curve of genus $4$ using lemma \ref{affine} and  Corollary \ref{infinity}. If it does, we compute a list with its number of points over $\bF_2{i}$ for $i=1,2,3,4$ and include this list in the NNumbersH list and the pair $[p,q]$ in CurvesH list. 
\end{enumerate}

\section{How we computed all isomorphism classes: The trigonal
case}

In case we have a trigonal curve of genus 4 over a field $k$,  it is
a well-known special case of Noether-Enriques theorem that its
canonical image can be express as the intersection of a quadratic
surface and a cubic surface in $\bP^3$. Moreover, the rank of
quadratic surface (as a quadratic form) must be $\ge 3$: it is
either non singular (rank 4) or of rank 3. The following result is
well-known, but we write the proof for lack of a suitable reference.

\begin{lem}\label{quadrics} Let $k$ be a \textbf{finite} field of characteristic 2, and
let $q\in k[x,y,z,t]$ an homogeneous polynomial of degree 2. Then,
after a suitable change of variables, there are only two
possibilities for $q$ it it has rank 4, and one if it has rank 3: 
$$q_1=xy+zt , \ q_2=xy+z^2+zt+at^2 \text{ and } q_3=xy+z^2$$
where $a\in k^*$ is a fixed element such that $x^2+x+a$ is irreducible, and being $q_3$ the singular case. 
\end{lem}
\begin{pf} It is a result by Arf that the
non-singular quadratic forms are determined by the so called Arf
invariant, which is an element in $L/\AS(L)$. Hence, over a finite field there are only two options. The Arf invariant of $q_1$ is $0$ and the one of $q_2$ is $a$, so they are not equivalent; this settles the non-singular case. Now, if $q\in k[x,y,z,t]$ is an homogeneous polynomial of degree 2 and rank $3$, then $q$ must be equivalent to either $q_3$ or to $ x^2+xy+ay^2+z^2 $ for some $a\in K^*$; but this two are equivalent by a change of variables of the form $z\mapsto x+by+z$, for some $b\in K^*$ with $b^2=a$. 
\end{pf}

We implemented the following method to list all non-hyperelliptic
genus 4 curves over a finite field $k$. For a fixed quadric $q$ as
in the lemma \ref{quadrics}
\begin{enumerate}
\item We compute the group  $G_q$  of all the matrices in $\PGL_4(k)$ that fix the quadric $q$.
\item We fix an order of the monomials of degree 3 in four variables (there are 20 such monomials). For such an order, we have an identification of all homogeneous degree 3
polynomials with 4 variables with the projective space $\bP^{19}_k$ of
dimension $19=20-1$.
\item We compute all the degree 3 homogeneous polynomials that are multiple of $q$ (by multiplying $q$ by an homogeneous linear polynomial),
and the set $S_q$ of the corresponding points in $\bP^{19}_k$. The set of initial candidates will be the points in $\bP^{19}(k)$ not in $S_q$).
\item For a given candidate $w$,(i.e  in $\bP^{19}(k)$ but not in $S_q$), we compute the set of points $S_w$ obtained by applying the change
of variables to the associated cubic $f_w$ corresponding to the
matrices in $G_q$, and also adding an element of $S_q$. The set
$S_w$ correspond, therefore, to all the cubic polynomials that,
together with $q$, determine schemes which are linearly equivalent.
\item We take out $S_w$ from the set of initial candidates, and we
pass $w$ (or a minimal element in $S_w$ in some sense) to the next
step.
\item Finally, we determine if the scheme corresponds to a curve of genus 4
(i.e., it is a curve, and it is non singular). We add such a point
to the list of final candidates.
\item We compute a list with its number of points over $\bF_2{i}$ for $i=1,2,3,4$ and include this list in the NNumbersNH list and the pair $[p,q]$ in CurvesNH list. 
\item We repeat until we exhaust the set of initial candidates.
\end{enumerate}

After we computed all the hyperelliptic and trigonal curves, we made a list NNumbers of all the possible lists of N-numbers ordered lexicographically, and we produce two lists with the same lenght than the NNumbers list containing a list of the hyperelliptic curves and of the trigonal curves respectively that in the index i have the same N-numbers than NNumbers[i].

%\newpage

\section{Tables}

In this section we collect some of the data corresponding to the curves. This data and the code in Magma \cite{magma} of the computations can be found in the github account of the author

\url{https://github.com/XavierXarles/Censusforgenus4curvesoverF2}

The following theorem collects some of the data we obtained after the computation. 

\begin{thm} There are 264 hyperelliptic curves of genus 4 over $\bF_2$ pairwise non-isomorphic. 

There are 780 trigonal curves of genus 4 over $\bF_2$ pairwise non-isomorphic. 

So, in total there are 1044 isomorphic classes of genus 4 curves over $\bF_2$ 

There are 620 classes modulo isogeny for jacobians of genus 4 curves. 

\end{thm}

The following table we give the number of curves up-to-isomorphy having a fixed number of points over $\bF_2$.

\begin{table}[h!]
\centering
\begin{tabular}{|l|l|l|l|l|l|l|l|l|l|l|}
\hline \textbf{Number of Points} & 0 & 1 & 2 &3 &4 &5 &6 & 7 &8 & \textbf{Total}\\
\hline \textbf{\# of hyperel. curves}& 9& 32& 58& 66& 58& 32& 9& 0& 0& 264 \\
\hline \textbf{\# of trigonal curves}& 31& 119& 202& 201& 117& 68& 30& 11& 1& 780  \\
\hline \textbf{Total number of curves}& 40& 151& 260& 267& 175& 100&
39& 11& 1 & 1044
\\ \hline
 \end{tabular}
  \caption{Curves for a given number of points.}
\label{table:a1}
\end{table}

For a given curve $C$ over $\bF_2$, we denote by $a_n(C)$ the number of points of degree exactly $i$; so, 
$$na_n(C)=|C(\bF_{2^n})\setminus \bigcup_{d\text{ divides } n, d<n} C(\bF_{2^d}) |$$
(since any point in $C(\bF_{2^n})$ not in an smaller subfield comes with $n$-conjugates). Note that $a_1=|C(\bF_2)|$. In the following tables we given the number of curves for the given $a_n$.

\begin{table}[h!]
\centering
\begin{tabular}{|l|l|l|l|l|l|l|l|l|l|}
\hline \textbf{$a_2$} & 0 & 1 & 2 &3 &4 &5 &6 & 7 & \textbf{Total}\\
\hline \textbf{\# of hyperel. curves}& 33& 61& 79& 61& 27& 3& 0& 0& 264 \\
\hline \textbf{\# of trigonal curves}& 78& 169& 211& 175& 103& 35& 8& 1& 780  \\
\hline \textbf{\# of curves}& 111& 230& 290& 236& 130& 38&
8& 1 & 1044
\\ \hline
 \end{tabular}
 \caption{Curves for a given number of degree 2 points.}
\label{table:a2}
\end{table}

\begin{table}[h!]
\centering
\begin{tabular}{|l|l|l|l|l|l|l|l|l|l|l|}
\hline \textbf{$a_3$} & 0 & 1 & 2 &3 &4 &5 &6 & 7 &8 & \textbf{Total}\\
\hline \textbf{\# of hyperel. curves}& 55& 28& 98& 28& 55& 0& 0& 0& 0 & 264 \\
\hline \textbf{\# of trigonal curves}&  53& 135& 195& 180& 109& 62& 36& 8& 2 & 780  \\
\hline \textbf{\# of curves}& 108& 163& 293& 208& 164& 62&
36& 8& 2 & 1044
\\ \hline
 \end{tabular}
 \caption{Curves for a given number of degree 3 points.}
\label{table:a3}
\end{table}

\begin{table}[h!]
\centering
\begin{tabular}{|l|l|l|l|l|l|l|l|l|l|l|l|l|}
\hline \textbf{$a_4$} & 0 & 1 & 2 &3 &4 &5 &6 & 7 &8 & 9 & 10 &\textbf{Total}\\
\hline \textbf{\# of hyp. c.}& 17& 20& 52& 39& 63& 38& 23& 12& 0& 0& 0& 264 \\
\hline \textbf{\# of trig. c.}& 19& 78& 135& 179& 152& 98& 68& 28& 16& 6& 1& 780  \\
\hline \textbf{\#of curves}& 36& 98& 187& 218& 215& 136&
91& 40& 16& 6& 1 & 1044
\\ \hline
 \end{tabular}
 \caption{Curves for a given number of degree 4 points.}
\label{table:a4}
\end{table}

It is well known that the numbers $(a_1,a_2,a_3,a_4)$ for genus $4$ curves determine the so called $L$-polynomial, which in turn determine the isogeny class of the Jacobian of the curve by Tate's theorem \cite{Tate}. But several curves can have the same $L$-polynomial (so isogeneous jacobians) (see next sections for some curious examples).

In the following table we collect how many curves (of the given type) share the same $L$-polynomial. For example, the first columns says that there are 403 $L$-polynomials of Jacobians of genus 4 curves that are not $L$-polynomials for hyperelliptic curves, and 99 that are not $L$-polynomials for trigonal curves; and the last zero of this column is tautological (since it says how many $L$-polynomials for jacobians of genus 4 curves are not $L$-polynomials for jacobians of genus 4 curves). On the contrary, the last 1 in the table says that there is one $L$-polynomial which is the $L$-polynomial of 7 non-isomorphic genus 4 curves (see the next section): it happens to be the same $L$-polynomial shared by 4 non-isomorphic hyperelliptic  genus 4 curves (which is the other 1 in the table).

\begin{table}[h!]
\centering
\begin{tabular}{|l|l|l|l|l|l|l|l|l|}
\hline \textbf{\# curves with given L} & 0 & 1 & 2 &3 &4 &5 &6 & 7\\
\hline \textbf{\# of hyperel. curves}& 403& 174& 40& 2& 1 & & & \\
\hline \textbf{\# of trigonal curves}& 99& 341& 128& 31& 15& 6 & &  \\
\hline \textbf{\#  of curves}& 0 & 361& 165& 49& 25& 15& 4&
1
\\ \hline
 \end{tabular}
  \caption{Number of curves for how many curves have the same $L$ polynomial. }
\label{table:L curves}
\end{table}

Note that the last line gives how many $L$-polynomials there are with $i$ curves in its isogeny class. So there are $361$ jacobians whose isogeny class contains only one element, and 259 jacobians (modulo isogeny) whose isogeny class contains more than one element; this means that for around 42\% of the L-polynomials, there is more than one curve with that polynomial. 

\section{Special Curves}

\begin{exmp} The minimal curve:
The curve which has the smallest number of degree $d$ points for all
$d\le 4$ is the trigonal curve given by the equations
$$ \begin{cases}
        X^2 + X~Y + Y^2 + Z~T =0\\
        Y^3 + X~Z^2 + Z^3 + X~Y~T + T^3=0
\end{cases}$$
It has $a_1=a_2=a_3=0$ and $a_4=1$.
\end{exmp}

\begin{exmp} The maximal curve: The only curve with 8 points over $\bF_2$ is
$$ \begin{cases}
   X~Y + Z~T =0 \\
        X~Y^2 + Y^3 + X^2~Z + Y^2~Z + X~Z^2 + X^2~T + Y^2~T + X~T^2
=0
\end{cases}
$$
It has, however, $a_2=a_3=0$, and $a_4=2$.
\end{exmp}

\begin{exmp} There is only one curve $C$ defined over $\bF_2$ that has the maximal number of degree 2 points, namely $a_2=7$. It reaches also the maximal number of points a genus 4 curve can have over $\bF_4$, namely $|C(\bF_4)|=15$, but having only one point over $\bF_2$; there are 7 other curves with this number of points, but they have more points over $\bF_2$.   It is the curve with equations
$$ \begin{cases}
 X^2 + X~Y + Y^2 + Z~T=0\\
                X^2~Y + X^2~T + X~Y^2 + X~T^2 + Z^3=0
\end{cases}
$$

%126
\end{exmp}

\begin{exmp} There is only one curve defined over $\bF_2$ that has the maximal number of degree 4 points, namely $a_4=10$. It reaches also the maximal number of points a genus 4 curve can have over $\bF_{2^4}$, namely $|C(\bF_{2^4})|=45$ (see \cite{Ho-La}) (there is another curve with that number of points, but with smaller $a_4$). 

It is the curve with equations
$$ \begin{cases}
 X^2 + X~Y + Y^2 + Z~T=0\\
                X^2~Y + X~Y^2 + X~Y~Z + X~Y~T + X~T^2 + Y~Z^2=0
\end{cases}
$$

%297
\end{exmp}

\begin{exmp} There are two curves defined over $\bF_2$ that has the maximal number of degree 3 points, namely $a_3=8$, and which they also coincide with the curves having the maximal number of points a genus 4 curve can have over $\bF_{2^3}$, namely $|C(\bF_{2^3})|=25$ (see \cite{Savitt}). 

They are the curves with equations
$$ \begin{cases}
X~Y + T^2=0\\
 X^3 + X^2~Z + X~Y^2 + X~Y~Z + X~Y~T + X~Z~T + Y^3 + Y^2~T + Z^3=0
\end{cases}
$$
and
$$ \begin{cases}
 X~Y + Z~T=0\\
  X^3 + X^2~Z + X^2~T + X~T^2 + Y^3 + Y^2~Z + Y^2~T + Y~T^2 + Z^3=0
\end{cases}
$$

%81 and 102
\end{exmp}

\begin{exmp} There is only one curve defined over $\bF_2$ that has the maximal number of degree 5 points, namely $a_5=14$, and which also coincides with the curve over $\bF_2$ having the maximal number of points over $\bF_{2^5}$, namely $|C(\bF_{2^5})|=71$ (so it has one point over $\bF_2$); it is also the known maximal number of points for genus 4 curve over $ \bF_{2^5}$ (but it is indecided if there is one with 72 points, which necessarely will not be defined over $\bF_{2}$, see \cite{Ho-La2}). 

It is the curve with equations
$$ \begin{cases}
 X^2 + X~Y + Y^2 + Z~T=0\\
        X^3 + X^2~Y + X~Y^2 + X^2~Z + Y^2~Z + Y~Z^2 + Z^3 + X^2~T + X~T^2=0
\end{cases}
$$

%85
\end{exmp}

There are also 4 curves having $a_5=0$, so the same number of points over $\bF_2$ and over  $\bF_{2^5}$. Two of them have isogeneous jacobian.

\begin{exmp}

There are 4 genus 4 hyperelliptic curves , pairwise non isomorphic,
with the same Weil Polynomial, so with isogenous Jacobian. They are
given by the following equations:
$$
\begin{cases}
y^2+ ( x^4 + x^3 + 1)y= x^9 + 1 \\
y^2+( x^4 + x^3 + 1)y= x^9 + x^8 + x\\
y^2+( x^4 + x^3 + 1)y=  x^9 + x^8 + x^3\\
y^2+( x^4 + x^3 + 1)y= x^9 + x^3 + x + 1 
\end{cases}
$$
That they share the same class modulo isogeny can be shown easily by counting points. Since $G_q=1$ for $q=x^4 + x^3 + 1$, to curves with equations $y^2+q~y=p_1$ and $y^2+q~y=p_2$ can be isomorphic only if $p_2=p_1+r^2+qr$ for some $r$ polynomial of degree $\le 5$. Since all $p_i$ have degree 9 and $r$ has degree 8, we can restrict the search to polynomials $r$ of degree $\le 4$. From this description it is clear than if $p_1$ has a monomial of type $x^8$ it has also $p_2=p_1+r^2+qr$. A direct computation shows that the 4 equations give non isomorphic curves. 

In fact, there 3 other genus 4 curves, now trigonal, with the same
Weil polynomial, given by the canonical equations
$$
\begin{array}{l}
\begin{cases}        X~Y = Z~T \\
        X^3 + X~Y^2 + Y^3 + X^2~Z + X~Y~Z + X~Z^2 + X~T^2=0
\end{cases} \\
\begin{cases}
        X~Y = Z~T \\
        X^3 + X^2~Y + X~Y^2 + Y^3 + X^2~Z + X~Y~Z + Y^2~Z + \\  \hspace{ 150pt} X~Z^2 + X^2~T +
            X~Y~T + X~T^2=0
\end{cases} \\
\begin{cases}        X~Y = T^2\\
        X^2~Y + Y^3 + X~Y~Z + Z^3 + X^2~T + X~Y~T + Y^2~T + X~Z~T=0
\end{cases}
\end{array}
$$

These curves have $|C(\bF_2)|=3$, $|C(\bF_4)|=|C(\bF_8)|=9$ and $|C(\bF_{16})|=21$, so a numbers $(a_1,a_2,a_3,a_4)=(3,3,2,3)$. The jacobian is ordinary, and $\Jac(C)[2]\cong \bZ/30\bZ$.
\end{exmp}

In the cases descrived above one can see that two curves with distinc quadrics $q_i$ in the trigonal case can have the same jacobian modulo isogeny. Even more, there are also (two) cases of curves with isogenous jacobian and all
the possible quadrics in the trigonal case and also containing an hyperelliptic curve. 

\begin{exmp}
For example, the jacobian of the hyperelliptic curve 
$$y^2+xy= x^9 +1 $$
is isogenous to the jacobian of 4 distinc trigonal curves, and for each of the quadrics $q_i$, $i=1,2,3$, there is one curve with that quadric part.
\end{exmp}


\begin{thebibliography}{99999}

\bibitem{LMFDB} The LMFDB Collaboration, The L-functions and Modular Forms Database, http://www.lmfdb.org, 2019.

\bibitem{Ho-La}  E. W. Howe and K. E. Lauter,
Improved upper bounds for the number of points on curves over finite fields,
Ann. Inst. Fourier (Grenoble) 53 (2003) 1677-1737.

\bibitem{Ho-La2}  	 E. W. Howe and K. E. Lauter,
New methods for bounding the number of points on curves over finite fields,
Geometry and Arithmetic (C. Faber, G. Farkas, and R. de Jong, eds.), European Mathematical Society, 2012, pp. 173-212

\bibitem{Liu} Q. Liu, Algebraic Geometry and Arithmetic Curves, Oxford University Press, 2002. 

\bibitem{magma} W. Bosma, J. J. Cannon, C. Fieker, A. Steel (eds.), Handbook of Magma functions, Edition 2.25 (2020), 5017 pages.

\bibitem{NaSa} E. Nart, D. Satoril, Hyperelliptic curves of genus
three over finite fields of even characteristic, Finite Fields Appl.
10 (2004) 198-220
\bibitem{Savitt} D. Savitt, 
The maximum number of points on a curve of genus 4 over $F_8$ is 25,
Canad. J. Math. 55 (2003) 331-352.
\bibitem{Tate} J. Tate. Endomorphisms of Abelian Varieties over Finite Fields. Inventiones Math. 2, 1966.
p. 134-144.
\end{thebibliography}
\end{document}